\documentclass{amsart}
\usepackage{amsmath,amsthm}
\usepackage{amsfonts,amssymb}

\newtheorem{thm}{Theorem}[section]

\newtheorem{lem}[thm]{Lemma}
\newtheorem{prop}[thm]{Proposition}

\newtheorem{defn}[thm]{Definition}

\theoremstyle{remark}

 \def\CA{{\mathcal A}}
 
 \def\CD{{\mathcal D}}
 
 \def\CH{{\mathcal H}}
 
 \def\CJ{{\mathcal J}}

 \def\CP{{\mathcal P}}     
 \def\CS{{\mathcal S}}

 \def\NN{{\mathbb N}}

 \def\RR{{\mathbb R}}
 \def\ZZ{{\mathbb Z}}

        \def\sign{\operatorname{sign}}
        \def\rad{\operatorname{rad}}

\newcommand{\wt}{\widetilde}
\newcommand{\wh}{\widehat}

\begin{document}

\title
 {Riesz transform and Riesz potentials for Dunkl transform}
 
\author{Sundaram Thangavelu and Yuan Xu}
\address{Stat-Math Division\\Indian Statistical Institute\\
8th Mile, Mysore Road\\Bangalore-560 059\\India.}
\email{veluma@isibang.ac.in}
\address{Department of Mathematics\\ University of Oregon\\
    Eugene, Oregon 97403-1222.}\email{yuan@math.uoregon.edu}

\date{\today}
\thanks{The work of YX was supported in part by the National Science 
Foundation under Grant DMS-0201669}

\begin{abstract}
Analogous of Riesz potentials and Riesz transforms are defined and studied
for the Dunkl transform associated with a family of weighted functions 
that are invariant under a reflection group. The $L^p$ boundedness of 
these operators is established in certain cases. 
\end{abstract}

\maketitle

\section{Introduction}
\setcounter{equation}{0}

For a family of weighted functions, $h_\kappa$, invariant under a finite 
reflection group, Dunkl transform is an extension of the Fourier transform
that defines an isometry of $L^2(\RR^d, h_\kappa^2)$ onto itself. The basic
properties of the Dunkl transforms have been studied by several authors,
see \cite{D92,Jeu,Ros,R03,So,TX,Tri} and the references therein. Giving the 
important role of Fourier transform in analysis, one naturally asks if it 
is possible to extend results established for the Fourier transform to the 
Dunkl transform. 

In analogous to the ordinary Fourier analysis, one can define a convolution
operator and study various summabilities of the inverse Dunkl transforms.
The convolution is defined through a generalized translation operator,
$\tau_y$, which plays the role of $f \mapsto f(\cdot - y)$ but is defined 
in the Dunkl transform side. The explicit expression of $\tau_y f$ is 
known only in some special cases and it is not a positive operator in 
general. In fact, even the $L^p$ boundedness of $\tau_y$ is not established 
in general. This is the main reason that only part of the results for the
Fourier transforms has been extended to the Dunkl transform at the moment. 

Recently, in \cite{TX}, the $L^p$ theory for convolution operators was 
studied. In particular, the $L^p$ boundedness of the convolution operator 
is established in the case that the kernel is a suitable radial function.
Furthermore, a maximal function is defined and shown to be of strong
type $(p,p)$ and weak type $(1,1)$. This provides a handy tool for 
extending some results from the Fourier transform to the Dunkl transform. 
In the present paper we study the analogous of the Riesz potentials and
the Riesz transforms for the Dunkl transform. We will study the boundedness
of the Riesz potentials as well as the related Bessel potentials. The 
Riesz transforms are examples of singular integrals. A general theory of 
singular integral for the Dunkl transform appears to be out of reach at 
the moment. We will prove 
the $L^p$ boundedness of the weighted Riesz transform only in a very 
special case of $d =1$ and $G = \ZZ_2$. Even in this simple case, however,
the proof turns out to be rather nontrivial.

The paper is organized as follows. In the next section we collect
the background materials. In Section 3 we recall the definition of the 
ordinary Riesz transforms and Riesz potentials, and prove a weighted 
$L^p$ boundedness for the Riesz potentials that will be used later 
in the paper. The weighted Riesz potentials and the Bessel potentials 
for the Dunkl transform will be studied in Section 4.  The weighted
Riesz transform is discussed in Section 5. 

Throughout this paper we use the convention that $c$ denotes a generic
constant, depending on $d$, $p$, $\kappa$ or other fixed parameters, 
its value may change from line to line.

\section{Preliminaries}
\setcounter{equation}{0}

\subsection{Dunkl Transform}

The Dunkl transform is associated to a weight function that is invariant 
under a reflection group. Let $G$ be a finite reflection group on $\RR^d$
with a fixed positive root system $R_+$, normalized so that 
$\langle v, v \rangle =2$ for all $v \in R_+$, where $\langle x,y \rangle$
denotes the usual Euclidean inner product. Let $\kappa$ be a nonnegative 
multiplicity function $v \mapsto \kappa_v$ defined on $R_+$ with the property 
that $\kappa_u = \kappa_v$ whenever $\sigma_u$ is conjugate to $\sigma_v$ 
in $G$; then $v \mapsto \kappa_v$ is a $G$-invariant function. The weight
function $h_\kappa$ is defined by 
\begin{equation}\label{eq:1.1}
h_\kappa(x) = \prod_{v \in R_+} |\langle x, v\rangle|^{\kappa_v}, \qquad
   x \in \RR^d.
\end{equation}
This is a positive homogeneous function of degree $\gamma_\kappa:=
\sum_{v \in R_+} \kappa_v$, and it is invariant under the reflection group 
$G$. 

To define the Dunkl transform we will also need the intertwining operator
$V_\kappa$. Let $\CD_j$ denote Dunkl's differential-difference operators 
defined by \cite{D89}
$$
  \CD_j f(x) = \partial_j f(x) + \sum_{v \in R_+} k_v
    \frac{f(x) -  f(x \sigma_v)} {\langle x, v\rangle}
        \langle v,\varepsilon_j\rangle, \qquad 1 \le j \le d,
$$
where $\varepsilon_1, \ldots, \varepsilon_d$ are the standard unit vectors of
$\RR^d$ and $\sigma_v$ denote the reflection with respect to the hyperplane 
perpendicular to $v$, $x \sigma_v:= x-2(\langle x,v \rangle /\|v\|^2) v$, 
$x \in \RR^d$. The operators $\CD_j$, $1\le j\le d$, map $\CP_n^d$ to 
$\CP_{n-1}^d$, where $\CP_n^d$ denotes the space of homogeneous polynomials 
of degree $n$ in $d$ variables, and they mutually commute; that is, 
$\CD_i \CD_j =\CD_j \CD_i$, $1 \le i,j \le d$. The intertwining operator 
$V_\kappa$ is a linear operator determined uniquely by 
$$
V_\kappa \CP_n \subset \CP_n, \quad V_\kappa1=1, \quad
\CD_i V_\kappa = V_\kappa \partial_i, \quad 1\le i\le d.
$$
The explicit formula of $V_\kappa$ is not known in general. For the group
$G = \ZZ_2^d$, $h_\kappa(x) = \prod_{i=1}^d |x_i|^{\kappa_i}$, it is an 
integral transform 
\begin{equation}\label{eq:1.2}
  V_\kappa f(x) = b_\kappa \int_{[-1,1]^d} f(x_1 t_1,\ldots, x_d t_d)
   \prod_{i=1}^d (1+t_i)(1-t_i^2)^{\kappa_i-1} dt.
\end{equation}
It is known that $V_\kappa$ is a positive operator \cite{Ros}; that is,
$p \ge 0$ implies $V_\kappa p\ge 0$.

Let $E(x,iy) = V_\kappa^{(x)} \left[ e^{i \langle x,y\rangle} \right]$,
$x ,y \in \RR^d$, where the superscript means that $V_\kappa$ is applied to 
the $x$ variable. For $f \in L^1(\RR^d, h_\kappa^2)$, the Dunkl transform is 
defined by 
\begin{equation} \label{eq:1.3}
\wh f(y) =  c_h \int_{\RR^d} f(x) E(x,-iy) h_\kappa^2(x)dx
\end{equation}
where $c_h$ is the constant defined by $c_h^{-1} = \int_{\RR^d} h_\kappa^2(x)
e^{-\|x\|^2/2}dx$. If $\kappa = 0$ then $V_\kappa = id$ and the Dunkl 
transform coincides with the usual Fourier transform. If $d =1$ and $G=\ZZ_2$,
then the Dunkl transform is related closely to the Hankel transform on the 
real line. 

Some of the properties of the Dunkl transform is collected below 
(\cite{D92,Jeu}).

\begin{prop} \label{prop:1.2}

\begin{enumerate}
\item For $f\in L^1(\RR^d, h_\kappa^2)$, $\wh f$ is in $C_0(\RR^d)$.

\item When both $ f $ and $\wh f $ are in $L^1(\RR^d, h_\kappa^2) $
we have the inversion formula
$$
  f(x) = \int_{\RR^d} E(ix,y) \wh f(y) h_\kappa^2(y)dy.
$$

\item The Dunkl transform extends to an isometry of 
$L^2(\RR^d, h_\kappa^2)$.

\end{enumerate}
\end{prop}

\subsection{Generalized translation operator}

Let $y \in \RR^d$ be given. The generalized translation operator
$f \mapsto \tau_y f$ is defined on $ L^2(\RR^d, h_\kappa^2) $ by
the equation
\begin{equation}\label{eq:1.4}
  \wh{\tau_y f}(x) = E(y,-ix) \wh f(x), \qquad x \in \RR^d.
\end{equation}
It plays the role of the ordinary translation $\tau_y f = f(\cdot -y)$ of
$\RR^d$, since the Fourier transform of $\tau_y$ is given by
$\wh {\tau_y f}(x) = e^{-i \langle x, y \rangle} \wh f(x)$. 

The generalized translation operator has been studied in \cite{Ros,R03,TX,Tri}.
The definition gives $\tau_y f$ as an $L^2$ function. Let us define
$$ 
\CA_\kappa(\RR^d) = \{ f \in  L^1(\RR^d; h_\kappa^2): \wh{f} \in
  L^1(\RR^d; h_\kappa^2) \}. 
$$
Then \eqref{eq:1.4} holds pointwise. Note that $\CA_\kappa(\RR^d)$ is 
contained in the intersection of $L^1(\RR^d;h_\kappa^2)$ and $ L^\infty $ 
and hence is a subspace of $L^2(\RR^d; h_\kappa^2)$. The operator $\tau_y $ 
satisfies the following properties: 

\begin{prop} \label{prop:3.2}
Assume that $ f \in \CA_\kappa(\RR^d)$ and $ g \in  L^1(\RR^d; h_\kappa^2)$ 
is bounded. Then
\begin{enumerate}
\item
$\displaystyle{\int_{\RR^d} \tau_y f(\xi) g(\xi) h_\kappa^2(\xi) d\xi
     =  \int_{\RR^d} f(\xi) \tau_{-y} g(\xi) h_\kappa^2(\xi) d\xi.}
$
\medskip
\item $\tau_y f(x) = \tau_{-x} f(-y)$.
\end{enumerate} 
\end{prop} 

A formula of $\tau_y f$ is known, at the moment, only in two cases. One 
is in the case of $G = \ZZ_2$ and $h_\kappa(x) = |x|^\kappa$ on $\RR$
(\cite{R95}) 
\begin{align} \label{eq:1.5}
\tau_y f(x) =& \frac{1}{2}
      \int_{-1}^1 f\left(\sqrt{x^2+y^2-2 x y t} \right)
  \Big(1+\frac{x-y}{\sqrt{x^2+y^2 - 2 x y t}} \Big) \Phi_\kappa(t)dt \\
& + \frac{1}{2} \int_{-1}^1 f\left(-\sqrt{x^2+y^2-2 x y t }\right)
  \Big(1-\frac{x-y}{\sqrt{x^2+y^2 - 2 x y t}} \Big) \Phi_\kappa(t)dt, \notag
\end{align}
where $\Phi_\kappa(t) = b_\kappa (1+t)(1-t^2)^{\kappa-1}$, from which also 
follows a formula of $\tau_y f$ in the case of $G = \ZZ_2^d$. The explicit 
formula implies the boundendess of $\tau_y f$. 
Let $\| \cdot\|_{\kappa,p}$ denote the norm of $L^p(\RR^d, h_\kappa^2)$. 

\begin{prop} \label{prop:1.3}
Let $G = \ZZ_2^d$. For $f\in L^p(\RR^d,h_\kappa^2)$, $1 \le p \le \infty$,  
$$
    \|\tau_yf \|_{\kappa,p} \le c \|f\|_{\kappa,p}.
$$
\end{prop}

Another case where a formula for $\tau_y f$ is known is when $f$ are radial 
functions, $f(x) = f_0(\|x\|)$, and $G$ being any reflection group 
(\cite{R03}), 
\begin{align} \label{eq:1.6}
  \tau_y f(x) = V_\kappa \left[ f_0\left(\sqrt{\|x\|^2+\|y\|^2 -
    2 \|x\|\ \|y\| \langle x', \cdot \rangle}\right) \right](y'),
\end{align}
from which it follows that $\tau_y f(x) \ge 0$ for all $y \in \RR^d$ if 
$f(x) = f_0(\|x\|) \ge 0$. 

Several essential properties of $\tau_y f$ is established for $f$ being 
radial functions. This is collected in the following proposition (\cite{TX}).
Let $L^p_{\rad}(\RR^d, h_\kappa^2)$ stands for the subspace of radial functions
in $L^p(\RR^d,h_\kappa^2)$.

\begin{prop}  \label{prop:1.6}
\begin{enumerate}
\item For every $ f \in L^1_{\rad}(\RR^d; h_\kappa^2) $,
$$ 
\int_{\RR^d} \tau_yf(x) h_\kappa^2(x) dx = \int_{\RR^d} f(x) h_\kappa^2(x) dx.
$$
\item For $1 \le p \le 2$, $ \tau_y:L_{\rad}^p(\RR^d, h_\kappa^2) 
 \rightarrow L^p(\RR^d; h_\kappa^2) $ is a bounded operator.
\end{enumerate}
\end{prop}

The generalized translation $\tau_y$ also satisfies the following 
property \cite{TX,Tri}: If $f$ is supported in $\{x: \|x\| \le B\}$ then 
$\tau_y f$ is supported in $\{x: \|x\| \le B+\|y\| \}$.  An important
consequence of this property is as follows (\cite{TX}): If $f \in 
C_0^\infty(\RR^d)$ then for $1 \le p \le \infty$ 
$$
  \lim_{y \to 0} \|\tau_y f - f \|_{\kappa,p} = 0.
$$

\subsection{The generalized convolution and maximal function}

The generalized translation operator can be used to define a convolution. 
For $f,g$ in $L^2(\RR^d,h_\kappa^2)$, we define
\begin{equation}\label{eq:convo} 
    (f *_\kappa g)(x) = \int_{\RR^d} f(y) \tau_x g^\vee(y) h_\kappa^2(y) dy
\end{equation}
where $g^\vee(y) = g(-y)$. 
Since $\tau_{x} g^\vee \in L^2(\RR^d,h_\kappa^2)$ the convolution is well
defined. 

This convolution has been considered by several authors (\cite{R03,TX,Tri}
and the refrences therein). It satisfies the basic properties 
$\wh {f *_\kappa g} = \wh f \cdot \wh g$ and $f *_\kappa g = g *_\kappa f$. 
Furthermore, if the generalized translation operator is bounded in norm, 
then the usual Young's inequality holds. For the general reflection group,
the following result is proved in \cite{TX}. 

\begin{thm} 
Let $ g $ be a bounded radial function in $ L^1(\RR^d;h_\kappa^2).$ Then
$f*_\kappa g$ initially defined in \eqref{eq:convo} on the intersection of 
$ L^1(\RR^d;h_\kappa^2)$ and $ L^2(\RR^d;h_\kappa^2) $ extends to all
$ L^p(\RR^d;h_\kappa^2)$, $1 \leq p \leq \infty $ as a bounded operator.
In particular,
\begin{equation}\label{eq:convo2}
   \|f*_\kappa g\|_{\kappa,p} \leq \|g\|_{\kappa,1} \|f\|_{\kappa,p}.
\end{equation}
\end{thm}

For $ f \in L^2 $ the maximal function $M_\kappa f$ is defined in \cite{TX} by
$$ 
  M_\kappa f(x) =  \sup_{r>0} \frac{1}{ d_\kappa r^{d+2\gamma_\kappa}} 
     |f *_\kappa \chi_{B_r}(x)|,
$$
where $ \chi_{B_r}$ is the characteristic function of the ball $B_r$
of radius $r$ centered at $0$ and 
$$
(d_\kappa)^{-1} = \int_{B_1} h_\kappa^2(y) dy =
(d+2\gamma_\kappa)\int_{S^{d-1}} h_\kappa^2(x) d\omega.
$$
The maximal function can also be written as
$$
  M_\kappa f(x) = \sup_{r>0} \frac{\left |
      \int_{B_r} \tau_y f(x) h_\kappa^2(y) dy \right| }
     {\int_{B_r} h_\kappa^2(y) dy}.
$$
Since $\tau_y \chi_{B_r}(x) \ge 0$ we have $M_\kappa f(x) \le M_\kappa |f|(x)$.
Furthermore, the following theorem holds: 

\begin{thm} \label{thm:1.4}
The maximal function is bounded on $L^p(\RR^d, h_\kappa^2)$ for $1 < p
\le \infty$; moreover it is of weak type $(1,1)$, that is, for $f \in
L^1(\RR^d,h_\kappa^2)$ and $a >0$,
$$
  \int_{E(a)} h_\kappa^2(x) dx \le \frac{c}{a} \|f\|_{\kappa,1}
$$
where $E(a) = \{ x: M_\kappa f(x) > a\}$ and $c$ is a constant independent
of $a$ and $f$.
\end{thm}

For $\phi \in L^1(\RR^d, h_\kappa^2)$ and $\varepsilon > 0$, the
dilation $\phi_\varepsilon$ is defined by
\begin{equation*}
 \phi_\varepsilon(x) = \varepsilon^{- (2\gamma_\kappa+d)}
    \phi(x/\varepsilon).
\end{equation*}
A change of variables shows that
$$
  \int_{\RR^d} \phi_\varepsilon (x) h_\kappa^2(x) dx =
      \int_{\RR^d} \phi (x) h_\kappa^2(x) dx, \qquad
      \hbox{for all $\varepsilon > 0$}.
$$
For convolution with an radial kernel, the following result is established
in \cite{TX}.

\begin{thm} \label{thm:6.2}
Let $ \phi \in \CA_\kappa(\RR^d) $ be a real valued radial function which
satisfies $ |\phi(x)| \le c (1+\|x\|)^{-d- 2\gamma_\kappa -1}.$ Then 
$$ 
\sup_{\varepsilon >0}|f*_\kappa \phi_\varepsilon(x)| \le c M_\kappa f(x).
$$
Consequently, $ f*_\kappa \phi_\varepsilon (x) \rightarrow f(x) $ for 
almost every $ x $ as $ \varepsilon $ goes to $ 0$ for all $ f $ in 
$ L^p(\RR^d;h_\kappa^2)$, $1 \leq p < \infty$.
\end{thm}

If $\tau_y$ is bounded in $L^p(\RR^d;h_\kappa^2)$ then the conditions in the
above theorem can be relaxed. At the moment this holds in the case of 
$G = \ZZ_2^d$ (\cite{TX}). 

\begin{thm} \label{thm:1.5}
Set $G = \ZZ_2^d$. Let $\phi(x) = \phi_0(\|x\|)  \in L^1(\RR^d,h_\kappa^2)$
be a radial function. Assume  that $\phi_0$ is differentiable,
$\lim_{r \to \infty} \phi_0(r) = 0$ and
$\int_0^\infty r^{2 \gamma_\kappa+ d} |\phi_0(r)| dr < \infty$, then
$$
  |(f *_\kappa \phi)(x)| \le c M_\kappa f(x) 
       \int_0^\infty r^{2 \gamma_\kappa+d} |\phi_0(r)| dr < \infty.
$$
In particular, if $\wh \phi(x) = \Phi(x)$ satisfies $\phi \in
L^1(\RR^d,h_\kappa^2)$ and $\Phi(0)=1$, then 
\begin{enumerate}
\item For $1 \le  p \le \infty$, $f *_\kappa \phi_\varepsilon$
converges to $f$ as $\varepsilon \to 0$ in $L^p(\RR^d,h_\kappa^2)$;
\item For $f \in L^1(\RR^d,h_\kappa^2)$,
$(f *_\kappa \phi_\varepsilon)(x)$  converges to $f(x)$ as
$\varepsilon \to 0$ for almost all $x \in \RR^d$.
\end{enumerate}
\end{thm}

\section{Ordinary Riesz transforms and Riesz potentials}
\setcounter{equation}{0}

In this section the notation $\wh f$ denote the ordinary Fourier 
transform on $\RR^m$. 

We recall the classical definition of the Riesz transforms and Riesz 
potentials. The Riesz transform, $R_jf$ ($1 \le j \le m$), for the ordinary 
Fourier transform on $\RR^m$ is defined by
$$ 
 R_j f(x) = \lim_{\varepsilon \to 0} c_j \int_{\|y\|\ge \varepsilon}
  f(x-y) \frac{y_j}{\|y\|^{m+1} } d y, 
  \qquad c_j = \frac{2^{m/2}}{\sqrt{\pi}} \Gamma\left(\frac{m+1}{2}\right).  
$$
It is a multiplier operator in the sense that 
$$
\wh{R_j f} (x) = -i\frac{x_j}{\|x\|} \wh f(x) \qquad \hbox{in $L^2(\RR^m)$}
$$
with respect to the ordinary Fourier transform. It is well known that 
$R_j f$ is a bounded operator from $L^p(\RR^m)$ to $L^p(\RR^m)$ for 
$1 < p < \infty$. 

If $f$ is a radial function, $f(x) = f_0(\|x\|)$, then the 
spherical-polar coordinates and the Funk-Hecke formula give us 
\begin{align}\label{eq:2.0}
R_j f(x) &= c \int_{\RR^m} f(x - y) \frac{y_j}{\|y\|^{ m +1}} dy = 
  c \int_0^\infty \int_{S^{d-1}} f_0(\|x - sy'\|) 
   y_j' d\omega(y') \frac{d s}{s} \\
 & = c x_j' \int_0^\infty \int_{-1}^1 f_0\left(
 \sqrt{\|x\|^2+ s^2 - 2\|x\| s t}\right) t (1-t^2)^{\frac{m-3}{2}} 
   dt   \frac{ds}{s}.  \notag
\end{align}
In particular, the Riesz transform of a radial function is also radial.
We introduce the following notation,  
$$
 \wt R_j f_0(\|x\|) = R_j f(x), \quad \hbox{when} \quad f(x) = f_0(\|x\|). 
$$
We will need the following proposition.

\begin{prop} \label{prop:3.1}
If $f(x) = f_0(\|x\|)$ and $f \in L^p(\RR^m)$ then 
$$
 \left(\int_0^\infty |\wt R_j f_0 (r)|^p r^{m-1} dr \right)^{1/p} \le
   c  \left(\int_0^\infty |f_0 (r)|^p r^{m-1} dr \right)^{1/p},
 \quad 1 < p < \infty.  
$$
\end{prop}

Evidently, using the spherical-polar coordinates $x = r x'$, this is a 
simple consequence of the boundedness of $R_j f$ on $L^p(\RR^m)$. 

The Riesz potential on $\RR^m$ is defined as the ordinary convolution of 
$f$ with the kernel $K_\alpha(x) = \gamma(\alpha)^{-1}\|x\|^{\alpha -m}$, 
$0 < \alpha < m$, 
\begin{equation} \label{eq:2.1}
f * K_\alpha (x) = \gamma(\alpha)^{-1} \int_{\RR^m} \|x - y\|^{-m +\alpha}
    f(y)dy,
\end{equation}
where $\gamma(\alpha)$ is a constant whose value we will not need.
It is well known (see, for example, \cite[p. 119]{S}) that for 
$1 < p < q < \infty$ where $1/q = 1/p - \alpha/m$,    
$$
   \|f * K_\alpha \|_q \le A_{p,q} \|f\|_p .
$$
The Riesz potential will be used in the Section 5. In particular,
we will need the following weighted inequality in the case of $\alpha =1$.

\begin{prop} \label{prop:2.2}
If $f(\|\cdot\|) \in L^p$, then 
$$
 \left( \int_{\RR^m} |f* K_1(x)|^p dx \right )^{1/p} 
   \le c  \left( \int_{\RR^m} |f(x)|^p \|x\|^p dx \right )^{1/p} 
$$
for all $p$ satisfying $p > m / (m-1)$, $m \ge 2$. 
\end{prop}

\begin{proof}
Let $M f$ be the Hardy-Littlewood maximal function defined by
$$
   M f(x) = \sup_{r > 0} \frac{1}{|B_r|} \int_{B_r} |f(x-y)| dy
$$
where $|B_r|$ denotes the volume of the ball $B_r$. We will make use
of the weighted norm ineqaulity
$$
   \int_{\RR^m} |M f(x)|^p w(x) dx \le c \int_{\RR^m} |f(x)|^p w(x) dx, 
$$
which holds whenever $w$ belongs to Muckenhoupt's $A_p$ class $A_p(\RR^m)$. 
It is known that $w(x) = \|x\|^\alpha \in A_p(\RR^m)$ whenever 
$- m < \alpha < m(p-1)$ (ses, for example, \cite[p. 218]{S95}).
In particular, $\|x\|^p \in A_p(\RR^d)$ for $p > m /(m-1)$ and 
$\|x\|^{- \delta p} \in A_p(\RR^d)$ for $0 < \delta < m/p$.

To prove the weighted ineqaulity we split $f*K_1$ as follows,
\begin{align*}
  (f * \|\cdot\|^{1-m})(x) & = 
     \int_{\|y\|< \|x\|} +
     \int_{\|x\|\le \|y\|\le 2 \|x\|}  +
     \int_{\|y\|> 2 \|x\|}  f(x- y) \|y\|^{1-m} dy \\
     & = T_1 f(x) + T_2 f(x) + T_3 f(x).           
\end{align*}
Evidently we have 
$$
  \|T_2 f(x)\| \le \int_{\|y\|\le 2 \|x\|} |f(x- y)| \,  \|y\|^{1-m} dy  
      \le c M f(x). 
$$
For $T_1$ we further split the integral to get 
\begin{align*}
|T_1f(x)|& \le \sum_{j=0}^\infty \int_{2^{-j-1} \|x\| \le \|y\|< 2^{-j} \|x\|} 
          | f(x- y)|\, \|y\|^{1-m} dy \\
 & \le \sum_{j=0}^\infty 2^{-j(1-m)} \|x\|^{1-m} 
        \int_{2^{-j-1} \|x\| \le \|y\|< 2^{-j} \|x\|} | f(x- y)| dy \\  
 & \le c \sum_{j=0}^\infty 2^{-j} \|x\| M f(x) \le c \|x\| M f(x). 
\end{align*}
Using the weighted ineqaulity of the maximal function it follows that, 
for $j=1,2$, 
$$
 \int_{\RR^d} |T_j f(x)|^p dx \le c \int_{\RR^d} |M f(x)|^p \|x\|^p dx 
    \le c \int_{\RR^d} |f(x)|^p \|x\|^p dx, 
$$  
as $\|x\|^p \in A_p(\RR^d)$ for $p > m/(m-1)$. 

To deal with $T_3f$ we choose $\delta$ such that $0 < \delta < m /p$ and
write the integral as
\begin{align*}
|T_3f(x)| & \le \sum_{j=1}^\infty \int_{2^j \|x\| \le \|y\|< 2^{j+1} \|x\|} 
          | f(x- y)| \, \|y\|^{1-m} dy \\
& = \sum_{j=1}^\infty \int_{2^j \|x\| \le \|y\|< 2^{j+1} \|x\|} 
          | f(x- y)| \, \|x-y\|^{1+\delta} \|x-y\|^{-1-\delta}\|y\|^{1-m} dy.  
\end{align*}
For $2^j \|x\| \le \|y\|< 2^{j+1} \|x\|$ we have $\|y\| \ge 2 \|x\|$ and
hence 
$$
 \|x-y\| \ge \|y\| - \|x\| = \|y\|/2 + (\|y\|/2 - \|x\|) \ge \|y\|/2 \ge
    2^{j-1}\|x\|  
$$ 
so that 
\begin{align*}
|T_3f(x)| & \le c\|x\|^{-\delta} \sum_{j=0}^\infty 2^{-j \delta}
   (2^j \|x\|)^{-m}
     \int_{\|x\| \le \|y\|< 2^{j+1} \|x\|} | f(x- y)| \|x-y\|^{1+\delta}dy \\  
 & \le c  \|x\|^{-\delta} \sum_{j=0}^\infty 2^{-j} M f_\delta(x) \\
 & \le c \|x\|^{-\delta} M f_\delta(x), 
\end{align*}
where $f_\delta(x)= f(x) \|x\|^{1+\delta}$. Then the weighted 
inequality of the maximal function implies, as $\|x\|^{-\delta p} \in 
A_p(\RR^d)$,
\begin{align*}
 \int_{\RR^d} |T_2 f(x)|^p dx & \le 
     c \int_{\RR^d} |M f_\delta(x)|^p \|x\|^{-\delta p}dx \\
   &  \le  c \int_{\RR^d} |f_\delta(x)|^p \|x\|^{-\delta p}dx  
      =  c \int_{\RR^d} |f(x)|^p \|x\|^p dx. 
\end{align*} 
This completes the proof. 
\end{proof} 

\section{Weighted Riesz Potentials and Bessel Potentials}
\setcounter{equation}{0}

In this section the notation $\wh f$ denotes the Dunkl transform of $f$.

\subsection{Riesz potentials} For $0 < \alpha < 2 \gamma_\kappa +d$,
the weighted Riesz potential, $I_\alpha^\kappa f$, is defined on $\CS$ by 
\begin{equation} \label{eq:RP} 
I_\alpha^\kappa f(x) = (d_\kappa^\alpha)^{-1} \int_{\RR^d} 
  \tau_y f(x)\frac{1}{\|y\|^{2\gamma_\kappa + d -\alpha}}h_\kappa^2(y) dy,   
\end{equation} 
where $d_\kappa^\alpha = 2^{-\gamma_\kappa-d/2 +\alpha}
 {\Gamma(\frac{\alpha}{2})}/{\Gamma(\gamma_\kappa+ \frac{d-\alpha}{2})}$. 

In order to derive the Dunkl transform of $I_\alpha^\kappa$, we start with 
a lemma, which is a little more general than what is needed. A homogeneous
polynomial $P$ is called an $h$-harmonics if $\Delta_h P =0$, where 
$\Delta_h = \CD_1^2+\ldots+\CD_d^2$ is called the $h$-Laplacian. Let 
$\CH_n^d(h_\kappa^2)$ denote the space of $h$-harmonics of degree $n$. It 
is known that 
$$
   \int_{S^{d-1}} P (x) q(x) h_\kappa^2(x) d\omega = 0
$$
whenever $P \in \CH_n^d(h_\kappa^2)$ and the degree $q$ is less than $n$.

\begin{lem} \label{lem:3.1}
For $P \in \CH_n^d(h_\kappa^2)$ and $0 < \Re \{\alpha\} < 2\gamma_\kappa+d$, 
the identity 
$$
\Big( \frac{P(x)}{ \|x\|^{2 \gamma_\kappa+d+n-\alpha}} \Big)\wh{ \Big .} 
      =  d_{n,\kappa} \frac{P(x)}{ \|x\|^{n+\alpha}}, \qquad  
   d_{n,\kappa}^\alpha = i^{-n} 2^{-\gamma_\kappa-d/2 +\alpha}
      \frac{\Gamma(\frac{n+\alpha}{2})} 
             {\Gamma(\gamma_\kappa+ \frac{n+d-\alpha}{2})}, 
$$
holds in the sense that
\begin{equation} \label{eq:3.0}
\int_{\RR^d} \frac{P(x)}{ \|x\|^{2 \gamma_\kappa+d+n-\alpha}}  
    \wh \phi(x) h_\kappa^2(x) dx
  = d_{n,\kappa}^\alpha \int_{\RR^d} \frac{P(x)}{\|x\|^{n+\alpha}}  
      \phi(x) h_\kappa^2(x) dx
\end{equation}
for every $\phi$ which is sufficiently rapidly decreasing at $\infty$,
and whose Dunkl transform has the same property. 
\end{lem}

\begin{proof}
If $P_n \in \CH_n^d(h_\kappa^2)$, then Theorem 5.7.5 of \cite{DX} shows
that 
$$
\left (P_n(x) e^{-\|x\|^2/2} \right) \wh{ \big .} = 
(-i)^n P_n(x) e^{-\|x\|^2/2}. 
$$
Since the Dunkl transform satisfies $\wh{f(s {\cdot})}(y) = 
s^{- 2\gamma_\kappa-d} \wh{f}(s^{-1} y)$, it follows that 
$$
  \left (P_n(x) e^{- s\|x\|^2/2} \right) \wh{ \big .} = 
(-i)^n s^{- n - \gamma_\kappa -d/2} P_n(x) e^{-\|x\|^2/(2s)}. 
$$
Let $\phi$ be a function that satisfies the property in the statement
of the lemma. For $s > 0$, the above formula leads to the relation 
\begin{align} \label{eq:3.1}
& \int_{\RR^d} P_n(x) e^{- s\|x\|^2/2} \wh \phi(x) h_\kappa^2(x) dx \\ 
& \qquad \qquad =  (-i)^n s^{- n - \gamma_\kappa -d/2} 
       \int_{\RR^d} P_n(x) e^{-\|x\|^2/(2s)} \phi(x)h_\kappa^2(x) dx. \notag
\end{align}
We then multiply the above equation by $s^{\beta -1}$, where 
$\beta = \gamma_\kappa + (d+n-\alpha)/2$, and integrate the result with 
respect to $s$ on $[0, \infty)$. Using 
$$
  \int_0^\infty s^{a -1} e^{-s \|x\|^2/2} ds = 2^a \Gamma(a)
    \|x\|^{- 2 a}
$$
and changing the order of the integrals, it is easy to see that 
this leads to
\begin{align*}
& 2^{\gamma_\kappa+(d+n-\alpha)/2} \Gamma(\gamma_\kappa+\tfrac{d+n-\alpha}{2})
\int_{\RR^d} \frac{P_n(x)}{\|x\|^{2 \gamma_\kappa + d+n-\alpha}} 
  \wh \phi(x) h_\kappa^2(x) dx \\
&\qquad\qquad\qquad =  (-i)^n 2^{(n+\alpha)/2} \Gamma(\tfrac{n+\alpha}{2}) 
 \int_{\RR^d} \frac{P_n(x)}{\|x\|^{n+\alpha}} \phi(x) h_\kappa^2(x) dx,
\end{align*}
which simplifies to the stated equation. In the above we can assume the
decay of $\phi$ and $\wh \phi$ is in the order of 
$$
 |\phi(x)| \le A(1+\|x\|)^{-d - 2 \gamma_\kappa} \quad \hbox{and}\quad
   |\wh \phi(x)| \le A(1+\|x\|)^{-d - 2 \gamma_\kappa}
$$
to ensure that the double integrals that occur above converge 
absolutely, so that the Fubini theorem applies. 
\end{proof}

\begin{prop}
Let $0 < \alpha < 2 \gamma_\kappa+d$. The identity 
\begin{equation}\label{eq:4.3}
  \wh I_\alpha^\kappa f (x) = \|x\|^{-\alpha} \wh f(x)
\end{equation}
holds in the sense that
$$ \int_{\RR^d} I_\alpha^\kappa f(x) g(x) h_\kappa^2 (x) dx 
 =  \int_{\RR^d} \wh f(x) \|x\|^{-\alpha} \wh g(x) 
    h_\kappa^2 (x) dx  
$$
whenever $f, g\in \CS$.  
\end{prop}

\begin{proof}
Setting $n = 0$ in Lemma \ref{eq:3.1} shows that the Dunkl transform
of $\|x\|^{-d-2\gamma_\kappa+\alpha}$ is the function $d_\kappa^\alpha 
\|x\|^{-\alpha}$ in the sense that 
$$ 
\int_{\RR^d} \|y\|^{-d - 2 \gamma_\kappa +\alpha} \phi (y) h_\kappa^2(y) dy
 =d_{\kappa}^\alpha \int_{\RR^d} \|y\|^{-\alpha} \wh \phi(y) h_\kappa^2(y)dy,
$$
where $\phi \in \CS$. Set $\phi(y) = \tau_y f(x)$ in the above identity
leads to 
$$ 
\int_{\RR^d} \|y\|^{-d - 2 \gamma_\kappa +\alpha}  
    \tau_y f(x) h_\kappa^2(y) dy
= d_{\kappa}^\alpha \int_{\RR^d} \|y\|^{-\alpha} E(-ix,y) \wh f(y) 
     h_\kappa^2(y) dy.
$$
Multiplying this identity by $g(x)$ and integrating, we obtain the stated
identity.
\end{proof}

Recall that $(-\Delta_h f)\,\wh{ }\, (x) = \|x\|^2 \wh f(x)$ (\cite{D92}) for
$f \in \CS$, the identity \eqref{eq:4.3} shows that the weighted Riesz 
potential can be defined as $( - \Delta_h)^{-\alpha/2} f$. The identity
\eqref{eq:4.3} also shows that 
$$
I_\alpha^\kappa ( I_\beta^\kappa f ) = I_{\alpha + \beta}^\kappa (f), \qquad
f \in \CS, \quad  \alpha,\beta > 0, \quad \alpha+ \beta < 2\gamma_\kappa+d.  
$$
$$
 \Delta_h (I_\alpha^\kappa f ) = I_\alpha^\kappa (\Delta_h f) 
    = - I_{\alpha-2}^\kappa(f), \qquad
   f \in \CS, \quad  2\gamma_\kappa + d > \alpha \ge 2,  
$$
which are extensions of familiar identities for the ordinary Riesz potentials. 
Next we consider the boundedness of $I_\alpha^\kappa$ as an operator from 
$L^q(\RR^d, h_\kappa^2)$ to $L^p(\RR^d, h_\kappa^2)$.  The following 
necessary condition holds:

\begin{prop}
If $\|I_\alpha^\kappa f\|_{\kappa,q} \le c \|f\|_{\kappa,p}$ for $f \in \CS$, 
then it is necessary that 
\begin{equation}\label{eq:4.4}
  \frac{1}{p} - \frac{1}{q} =  \frac{\alpha}{2\gamma_\kappa+d}. 
\end{equation}
\end{prop}

\begin{proof}
Let $f_s(x) = e^{-s^2 \|x\|^2}$. Using the fact 
$\tau_y f_s(x) = e^{-s^2(\|x\|^2+\|y\|^2)} E(2s^2x,y)$, 
a change of variable shows that 
\begin{align*}
 I_\alpha^\kappa (e^{-s^2 \|x\|^2})(x) = s^{-\alpha} 
       I_\alpha^\kappa (e^{- \|x\|^2})(s x). 
\end{align*} 
Consequently, setting $f(x) = 3^{-\|x\|^2}$, a changing of variables shows 
that 
$$
\|I_\alpha^\kappa f_s\|_{\kappa,q} = s^{-\alpha - (2 \gamma_\kappa +d)/q} 
  \|I_\alpha^\kappa f\|_{\kappa,q} \quad \hbox{and} \quad 
 \|f_s\|_{\kappa,p} = s^{- (2 \gamma_\kappa +d)/p} \|f\|_{\kappa,p} 
$$
for all $s > 0$. Considering $s \to 0$ and $s \to \infty$ shows that 
if $\|I_\alpha f\|_{\kappa,q} \le c \|f\|_{\kappa,p}$ holds, then we must 
have $ \alpha + (2\gamma_\kappa+d)/q -(2\gamma_\kappa+d)/p=0$, which gives
\eqref{eq:4.4}.
\end{proof}

The main result on the weighted Riesz potential is the following theorem.

\begin{thm} \label{thm:4.4}
Let $G = \ZZ_2^d$. Let $\alpha$ be a real number such that $0 < \alpha < 
2 \gamma_\kappa+d$ and let $1 \le p < q < \infty$ satisfies \eqref{eq:4.4}. 
\begin{enumerate}
\item For $f\in  L^p(\RR^d, h_\kappa^2)$, $p >1$, 
$$
  \|I_\alpha^\kappa f \|_{\kappa,q} \le c \|f\|_{\kappa,p}.
$$ 
\item For $f \in L^1(\RR^d, h_\kappa^2)$, the mapping 
$f \mapsto I_\alpha^\kappa f$ is of weak type $(1,q)$; that is,
$$
  \int_{\{x:|I_\alpha^\kappa(x)| > \sigma\}} h_\kappa^2(x) dx \le 
   c \left(\frac{\|f\|_{\kappa,1}}{\sigma} \right)^q.
$$
\end{enumerate}
\end{thm}

\begin{proof}
Let $R > 0$ be fixed. We write the operator as a sum of two terms, 
\begin{align*}
 I_\alpha^\kappa f(x) & = (d_\kappa^\alpha)^{-1} \int_{ \{x:\|x\|\le R\}} 
  \tau_y f(x)\frac{1}{\|y\|^{2\gamma_\kappa + d -\alpha}}h_\kappa^2(y) dy \\
 & + (d_\kappa^\alpha)^{-1} \int_{ \{x:\|x\| \ge R\}} 
  \tau_y f(x)\frac{1}{\|y\|^{2\gamma_\kappa + d -\alpha}}h_\kappa^2(y) d y  
  : = S_1 f(x) + S_2 f(x).
\end{align*} 
For $S_1 f$, we use the maximal function and the Theorem \ref{thm:1.5} 
to get the estimate
\begin{align*}
 |S_1 f(x)| \le c M_\kappa f(x) \int_0^R r^{2\gamma_\kappa + d} 
     \frac{d}{dr}\left( r^{-2\gamma_\kappa-d + \alpha}\right) dr
   \le  \frac{c}{\alpha} R^\alpha M_\kappa f(x).
\end{align*}
To estimate $S_2 f$ we use Proposition \ref{prop:1.3} which states that
$\|\tau_y f\|_{\kappa,p} \le c \|f\|_{\kappa,p}$. Let $p' = p/(p-1)$. Then 
H\"older's inequality shows that 
\begin{align*}
  |S_2 f(x)|&\, \le \left(\int_{y:\|y\|\ge R} 
    \left(\|y\|^{-2\gamma_\kappa -d+\alpha}\right)^{p'}
    h_\kappa^2(y) dy \right)^{1/p'} \|\tau_{-x} f\|_{\kappa,p} \\
 & \le c R^{- (d+2\gamma_\kappa)/p + \alpha} \|f\|_{\kappa,p}
   = c R^{- (d+2\gamma_\kappa)/q} \|f\|_{\kappa,p},
\end{align*}
where the last step follows from \eqref{eq:4.4}. Together, the two 
estimates show that 
$$
 |I_\alpha^\kappa f(x)| \le c \left(R^\alpha M_\kappa f(x)
     + R^{- (d+2\gamma_\kappa)/q} \|f\|_{\kappa,p}\right) 
$$
for all $R > 0$. Choosing $R = (M_\kappa f(x) /
 \|f\|_{\kappa,p})^{-p/(d+2\gamma_\kappa)}$ and using \eqref{eq:4.4},
we obtain the inequality
$$
 |I_\alpha^\kappa f(x)| \le c \left(M_\kappa f(x)\right)^{p/q} 
     \|f\|_{\kappa,p}^{1-p/q}.
$$   
Consequently, for $p > 1$, we can use the $L^p$ boundedness of 
the maximal function in Theorem \ref{thm:1.4} to conclude that 
$$
  \|I_\alpha^\kappa f \|_{\kappa,q} \le 
   c  \|M_\kappa f(x)\|_{\kappa,p}^{p/q} \|f\|_{\kappa,p}^{1-p/q}
   \le c \|f\|_{\kappa,p}.  
$$

For $p =1$, we use the weak type $(1,1)$ inequality of the maximal
function to get
\begin{align*}
\int_{\{x: I_\alpha^\kappa(x)\ge \sigma\}} h_\kappa^2(x) dx & \,\le 
 \int_{\{x: c  \|M_\kappa f(x)\|_{\kappa,p}^{p/q} \|f\|_{\kappa,p}^{1-p/q}
 \ge \sigma \}} h_\kappa^2(x) dx \\ 
 & \le c \left(\frac{\|f\|_{\kappa,1}^{1-1/q}}{\sigma} \right)^q 
    \|f\|_{\kappa,1} =  c \left(\frac{\|f\|_{\kappa,1}}{\sigma} \right)^q. 
\end{align*}
The proof is completed.
\end{proof}

The proof shows that if $\tau_y$ is a bounded operator from 
$L^p(\RR^d, h_\kappa^2)$ to itself  for some other reflection group, then 
the conclusion of the theorem will hold for that group. At the moment, the 
theorem holds only if $G = \ZZ_2^d$ or if $f$ are radial functions and 
$1 \le p \le 2$ by Propostion \ref{prop:1.6}. 

\subsection{Weighted Bessel potentials}
The Bessel potentials is closely related to the Riesz potential. The kernel
functions for the Bessel potentials have essentially the same local behavior 
as that of the Riesz potentials as $\|x\| \to 0$, but have much better 
behavior for $\|x\|$ large. In analogous to the ordinary Fourier transform,
the weighted Bessel potentials, $\CJ_\alpha$, can be defined by
$$
 \CJ_\alpha^\kappa = (I - \Delta_h)^{-\alpha/2}, \qquad   \alpha > 0. 
$$

To be more precise, we define $\CJ_\alpha^\kappa$ as a convolution operator
$$
\CJ_\alpha^\kappa  f = f *_\kappa G_\alpha^\kappa, \quad \hbox{where}
  \quad   \wh G_\alpha^\kappa (x) = (1+\|x\|^2)^{-\alpha/2} 
$$
for $f \in \CS$. The following position gives an explicit expression for
$G_\alpha^\kappa$. 

\begin{prop}
Let $G_\alpha^\kappa$ be defined as above. Then $G_\alpha^\kappa(x) \ge 0$ 
for all $x \in \RR$, $G_\alpha^\kappa \in L^1(\RR^d, h_\kappa^2)$, and 
\begin{equation} \label{eq:G}
  G_\alpha^\kappa (x) =  \frac{1}{\Gamma(\alpha/2)}
      \int_0^\infty e^{-t} e^{-\|x\|^2/(4t)} t^{-\gamma_\kappa+(\alpha - d)/2}
       \frac{dt}{t}. 
\end{equation}
\end{prop}

\begin{proof}
This follows as in the ordinary Bessel potentials (\cite[p. 132]{S}).
Let us work backward and start with \eqref{eq:G}. Evidently then 
$G_\alpha^\kappa(x) > 0$. Furthermore, since 
$$
 c_h \int_{\RR^d}  e^{-\|x\|^2/(2t)} h_\kappa^2(x) dx = 
  (2t)^{\gamma_\kappa+d/2} c_h \int_{\RR^d}  e^{-\|u\|^2/2} h_\kappa^2(u) du = 
     t^{\gamma_\kappa+d/2}
$$
the Fubini's theorem applied to \eqref{eq:G}, which shows 
$$
 c_h \int_{\RR^d} G_\alpha^\kappa (x) h_\kappa^2(x) dx = 
  \frac{1}{\Gamma(\alpha/2)} \int_0^\infty t^{\alpha/2 -1} e^{-t} dt = 1.
$$
Using the fact that $\left(e^{-\|x\|^2 /(4t)}\right)\wh{ } = 
e^{-t \|x\|^2}$, it follows that 
$$
  \wh G_\alpha^\kappa(x) = \frac{1}{\Gamma(\alpha/2)}
    \int_0^\infty e^{-t} e^{- t \|x\|^2)} t^{\alpha/2}  \frac{dt}{t}
    = (1+\|x\|^2)^{-\alpha/2}, 
$$
where the interchange of integrals can be easily justified by Fubini's theorem.
\end{proof}

The behavior of $G_\alpha^\kappa$ is described in the following lemma.

\begin{lem} \label{lem:G}
For $\alpha > 0$, 
$$
G_\alpha^\kappa (x) \le c \left(1+\|x\|^{- 2\gamma_\kappa - d +\alpha}\right)
   e^{-\|x\|/2}, \qquad  \|x\| > 0.  
$$
\end{lem}

\begin{proof}
The elementary inequality $t+ r^2/(4t)\ge 2 \sqrt{t} \sqrt{r^2/4t} = r$ 
leads to 
\begin{align*}
 G_\alpha^\kappa(x)
  \le \frac{1}{\Gamma(\alpha/2)} e^{-\|x\|/2} 
     \int_0^\infty e^{-\frac{1}{2} (t + \frac{\|x\|^2}{4t})}
         t^{-\gamma_\kappa+(\alpha - d)/2}  \frac{dt}{t}. 
\end{align*}
To estimate the integral we split it into two parts. Changing variable gives
\begin{align*}
&  \int_0^{\|x\|^2} e^{-\frac{1}{2} (t + \frac{\|x\|^2}{4t})}
         t^{-\gamma_\kappa+(\alpha - d)/2}  \frac{dt}{t} 
   \le \int_0^{\|x\|^2} e^{- \frac{\|x\|^2}{8t}}
         t^{-\gamma_\kappa+\frac{\alpha - d}{2}}  \frac{dt}{t}\\ 
& \qquad =\|x\|^{-2\gamma_\kappa + \alpha -d} 
      \int_0^1 e^{- \frac{1}{8u}}
         u^{-\gamma_\kappa+ \frac{\alpha - d}{2}}  \frac{du}{u}  
    \le c \|x\|^{-2\gamma_\kappa + \alpha -d}.  
\end{align*}
Furthermore, if $\|x\| \ge 1$, then 
\begin{align*}
  \int_{\|x\|^2}^\infty e^{-\frac{1}{2} (t + \frac{\|x\|^2}{4t})}
         t^{-\gamma_\kappa+(\alpha - d)/2}  \frac{dt}{t} 
   \le \int_{\|x\|^2}^\infty e^{- \frac{t}{2}}
         t^{-\gamma_\kappa+ \frac{\alpha - d}{2}}  \frac{dt}{t}\le c 
\end{align*}
and if $\|x\| \le 1$, then
\begin{align*}
&  \int_{\|x\|^2}^\infty e^{-\frac{1}{2} (t + \frac{\|x\|^2}{4t})}
         t^{-\gamma_\kappa+(\alpha - d)/2}  \frac{dt}{t} 
   \le \int_{\|x\|^2}^\infty e^{ -\frac{t}{2} }
         t^{-\gamma_\kappa+(\alpha - d)/2}  \frac{dt}{t} \\
&   \le c \left (1 + \int_{\|x\|^2}^1  t^{-\gamma_\kappa+(\alpha - d)/2}
    \frac{dt}{t} \right)
 \le c \left (1 + \|x\|^{-2\gamma_\kappa + \alpha -d} \right).  
\end{align*}
Putting these estimates together proves the stated inequality.
\end{proof}

This shows, in particular, that $G_\alpha^\kappa$ behaviors as 
$\|x\|^{- 2\gamma_\kappa -d + \alpha}$ for $\|x\| \to 0$, same as 
the kernel for the Riesz potentials. 

\begin{thm}
Under the same assumption, the conclusion of Theorem \ref{thm:4.4} holds for 
Bessel potentions. 
\end{thm}

The proof is essentially the same. We will not repeat it. In stead, we state
the following theorem which holds for all reflection groups.

\begin{thm}
Let $\alpha > 0$. 
\begin{enumerate}
\item The Bessel potentials are bounded operators from $L^p(\RR^d,h_\kappa^2)$
to itself for $1 \le p \le \infty$. 
\item For $f \in L^1(\RR^d, h_\kappa^2)$, 
$$
  | \CJ_\alpha^\kappa f(x) | \le c M_\kappa f(x), \qquad x \in \RR^d. 
$$ 
\end{enumerate}
\end{thm}

\begin{proof}
Since $G_\alpha^\kappa(x)$ is a radial function, let us write 
$G_\alpha^\kappa(r)$ for the function defined on $\RR_+$. 
The estimate in Lemma \ref{lem:G} shows that $G_\alpha^\kappa \in 
L^1(\RR^d, h_\kappa^2)$ and it has integral $1$, so that the first 
part of the theorem follows from Theorem \ref{thm:1.5}. 

For $\alpha > 0$, the estimate in Lemma \ref{lem:G} shows that the 
conditions on $\phi$ of Theorem \ref{thm:6.2} is satisfied which gives 
the second part. 
\end{proof}

Such a theorem does not hold for the Riesz potentials. The definition of the 
Bessel potentials also shows that $\CJ_\alpha^\kappa$ satisfies
$$
 \CJ_\alpha^\kappa \cdot  \CJ_\beta^\kappa =  \CJ_{\alpha+\beta}^\kappa, 
      \qquad \alpha > 0, \quad \beta >0.
$$

\section{Weighted Riesz Transforms}
\setcounter{equation}{0}

\subsection{Definition of Riesz transforms}
For $P \in \CH_n^d(h_\kappa^2)$, we consider the transform 
$$
 T^\kappa f(x) = \lim_{\varepsilon \to 0}  \int_{\|y\|\ge \varepsilon}
  \tau_y f(x) \frac{P(y)}{\|y\|^{\gamma_\kappa+d+n} } h_\kappa^2(y) d y  
$$
defined for $f \in L^2(\RR^d, h_\kappa^2)$. The equation \eqref{eq:3.1} 
and the Plancherel theorem shows that $T^\kappa f$ is well defined. 
We want to show that $T^\kappa$ is a multiplier operator under the Dunkl 
transform. A linear operator $T f$ is a multiplier operator if 
$\wh{T f}(x) = m(x) \wh{f}(x)$ in the sense that 
$$
   T f(x) = \int_{\RR^d}  m(y) \wh f(y) E(x,iy) h_\kappa^2(y) dy
$$
for $f$ sufficiently smooth and having compact support. 

\begin{thm}
Let $P \in \CH_n^d(h_\kappa^2)$, $n \ge 1$. Then the multiplier corresponding 
to the transform $T^\kappa f$ is given by $d_{n,\kappa} P(x) /\|x\|^k$, 
$d_{n,\kappa}= i^{-n} 
2^{- \gamma_\kappa-d/2} \Gamma(\frac{n}{2})/\Gamma(\gamma_\kappa+ 
\frac{n+d}{2})$. 
\end{thm}

\begin{proof}
Since $P_n \in \CH_n^d(h_\kappa^2)$, its integral with resepct to 
$h_\kappa^2$ on $S^{d-1}$ is zero. Hence,
$$
 \int_{\|x\| \le 1} \frac{P_n(x)}{\|x\|^{2\gamma_\kappa+d+n-\alpha}}
    h_\kappa^2(x) dx = \int_0^1 r^{\alpha -1} dr
       \int_{S^{d-1}} P_n(x') h_\kappa^2(x') d\omega(x') =0. 
$$
Let $\phi$ be a function that satisfies the condition in the Lemma 
\ref{lem:3.1} and 
the additional assumption that $\wh \phi$ is differentiable near 
origin. Then we can write
\begin{align*}
\int_{\RR^d} \frac{P_n(x)}{\|x\|^{2\gamma_\kappa+d+n-\alpha}} 
     \wh \phi(x) h_\kappa^2(x)dx 
 = &  \int_{\|x\| \ge 1} \frac{P_n(x)}{\|x\|^{2\gamma_\kappa+d+n-\alpha}}
     \wh \phi(x) h_\kappa^2(x)dx \\
 + & \int_{\|x\| \le 1} \frac{P_n(x)}{\|x\|^{2\gamma_\kappa+d+n-\alpha}}
   \left[\wh \phi(x) - \wh \phi(0)\right] h_\kappa^2(x)dx. 
\end{align*}
Since $[\wh \phi(x) - \wh \phi(0)]/\|x\|$ is locally integrable,
we can take limit $\alpha \to 0$ to get 
\begin{align*}
& \int_{\|x\| \le 1} \frac{P_n(x)}{\|x\|^{2\gamma_\kappa+d+n}}
  \left[\wh \phi(x) - \wh \phi(0)\right] h_\kappa^2(x)dx \\
&\qquad \qquad \qquad \qquad 
 = \lim_{\varepsilon \to 0}  \int_{\varepsilon \le \|x\| \le 1} 
  \frac{P_n(x)}{\|x\|^{2\gamma_\kappa+d+n}}\wh \phi(x)h_\kappa^2(x)dx. 
\end{align*}
Consequently, we conclude that 
\begin{equation}\label{eq:3.2}
 \lim_{\alpha \to 0+} \int_{\RR^d} \frac{P_n(x)}
 {\|x\|^{2\gamma_\kappa+d+n-\alpha}}
  \wh \phi(x) h_\kappa^2(x)dx
= \lim_{\varepsilon \to 0}  \int_{\|x\| \ge \varepsilon} 
  \frac{P_n(x)}{\|x\|^{2\gamma_\kappa+d+n}}\wh \phi(x)h_\kappa^2(x)dx. 
\end{equation}
Let $f$ be any sufficiently smooth function with compact support. For
a fixed $x$, set $\wh \phi(y) = \tau_y f(x) = \tau_{-x} f(-y)$. Then
\begin{align*} 
(\wh{\phi})\,\wh{\big.} \, (y) & =
  \int_{\RR^d} \tau_{-x} f(-z) E(z,-i y)h_\kappa^2(z) dz 
   = \int_{\RR^d} \tau_{-x} f(z) E(z, i y)h_\kappa^2(z) dz \\
  & = \wh {\tau_{-x} f} (-y) = E(-y,ix) \wh f(-y).
\end{align*}
Hence, it follows that 
$
 \phi(y) = (\wh{\phi})\, \wh{ {} }\, (- y) = E(y,ix) \wh f(y),
$
so that by \eqref{eq:3.0} and \eqref{eq:3.2}
$$
\lim_{\varepsilon \to 0}  \int_{\|x\| \ge \varepsilon} 
 \frac{P_n(y)}{\|y\|^{2\gamma_\kappa+d+n}} \tau_y f(x) h_\kappa^2(y)dy
  =  d_{n,\kappa}  \int_{\RR^d} 
   \frac{P_n(y)}{\|y\|^n} E(x,iy) \wh f(y) h_\kappa^2(y) dy. 
$$
By the definition of the multiplier $m$, 
we conclude that $m(y) = d_{n,\kappa} P_n(y) /\|y\|^n$.
\end{proof}

The proof of this theorem follows the argument for the ordinary Fourier 
transform as given in \cite[p. 73 -- 74]{S}.  

The special case that $P(x) = x_j$ defines the weighted Riesz transform.

\begin{defn} 
For $f \in L^2(\RR^d, h_\kappa^2)$ the Riesz transform $R_j^\kappa f$ 
is defined by
$$
 R_j^\kappa f(x) = \lim_{\varepsilon \to 0} c_j \int_{\|y\|\ge \varepsilon}
  \tau_y f(x) \frac{y_j}{\|y\|^{2 \gamma_\kappa+d+1} } h_\kappa^2(y) d y, 
$$
where $1 \le j \le d$ and $c_j = {2^{\gamma_\kappa + d/2}}
  {\Gamma(\gamma_\kappa +(d+1)/2)}/ {\sqrt{\pi}}$.    
\end{defn}

\begin{thm} The Riesz transform is a multiplier operator with 
$$
\wh{R_j^\kappa f}(x) = -i\frac{x_j}{\|x\|} \wh f(x), \qquad 1 \le j \le d.
$$
\end{thm}

\begin{proof}
Since $x_j \in \CH_1^d(h_\kappa^2)$, this is just the previous theorem
with $P(x) = x_j$.  
\end{proof}

\subsection{The boundedness of weighted Riesz transform}
The Riesz transforms are important singular integral operators. One would
like to prove the boundedness of the weighted Riesz transforms, just as in
the case of the ordinary Riesz transforms. This, however, turns out to be a 
rather difficult task. The effort is hindered by the lack of information on 
$\tau_y f$. Furthermore, currently no theory of singular integrals 
with reflection invariant weight functions is available.

For some special parameters, however, the weighted Riesz transforms on
the radial functions can be related to the classical Riesz transforms. 
This is given in the next proposition.

\begin{prop}
If $f(x) = f_0(\|x\|)$ is a radial function in $L^p(\RR^d, h_\kappa^2)$
and $2 \gamma_\kappa \in \NN_0$, then for $1 < p < \infty$,
$$
  \|R_j^\kappa f \|_{\kappa,p} \le c \|f\|_{\kappa,p}.
$$
\end{prop}

\begin{proof}
Since $f$ is radial, the explicit formula of $\tau_y f$ in \eqref{eq:1.5}
and the Funk-Hecke formula (\cite{X00}) shows that
\begin{align*}
R_j^\kappa f(x) &= c_h 
    \int_{\RR^d} \tau_y f(x) \frac{y_j}{\|y\|^{d+2\gamma_\kappa +1}}
     h_\kappa^2(y) dy \\
&  = c_h \int_0^\infty \int_{S^{d-1}} V_\kappa f_0\left(
   \sqrt{\|x\|^2+ s^2 - 2\|x\| s \langle x, \cdot \rangle }\right)
   (y') y_j' h_\kappa^2(y') d\omega(y') \frac{ds}{s}  \\
&  = c\, x_j' \int_0^\infty \int_{-1}^1 f_0\left(
 \sqrt{\|x\|^2+ s^2 - 2\|x\| s t}\right) t (1-t^2)^{\gamma_\kappa+(d-3)/2} 
   dt   \frac{ds}{s}.  
\end{align*}
Therefore, since $2 \gamma_\kappa \in \NN_0$, we conclude by \eqref{eq:2.0}
that $R_j^\kappa f(x) = c \wt R_j f_0( \|x\|)$, where $\wt R_j$ corresponds
to the ordinary Riesz transform $R_j f$ defined on $\RR^m$ with 
$m = d + 2 \gamma_\kappa$. Consequently, by Proposition \eqref{prop:3.1},
we have 
\begin{align*}
 \int_{\RR^d} |R_j^\kappa f(x)|^p h_\kappa^2(x)dx 
& =  c \int_0^\infty r^{2 \gamma_\kappa + d-1} |\wt R_j f_0(r)|^p dr \\ 
& =  c \int_0^\infty r^{m-1} |\wt R_j f_0(r)|^p dr  \\ 
& \le c \int_{\RR^m} |f_0(\|x\|)|^p dx = c \|f\|_{\kappa,p}^p
\end{align*}
which completes the proof.
\end{proof}

In the rest of this section, we consider only the case $d=1$ and $G= \ZZ_2$,
for which the weight function is simply $h_\kappa(x) = |x|^\kappa$ and the 
Riesz transform is the integral operator on the real line
$$
R^\kappa f(x) = 
 c_\kappa \int_{-\infty}^\infty \tau_y f(x) \frac{y}{\|y\|^{2 \kappa+2} } 
  h_\kappa^2(y) d y, \qquad c = \Gamma(\kappa +1) / \sqrt{\pi}, 
$$
where the integral holds in the principle value sense. There is a multiplier 
theorem in this setting of Dunkl transform on the real line (\cite{So}).
However, it does not apply to the Riesz transform. We have 
the following result. 

\begin{thm} 
Let $G = \ZZ_2$. If $f \in L^p(\RR^d, h_\kappa^2)$ and $2 \gamma_\kappa 
\in \NN_0$, then for $1 < p < \infty$,
$$
  \|R^\kappa f \|_{\kappa,p} \le c \|f\|_{\kappa,p}.
$$
\end{thm}

\begin{proof}
In this case $f$ is radial means that $f$ is even, so that the 
stated result holds for $f$ being even. Every function $f$ on
$\RR$ can be split as $f = f_e + f_o$ where $f_e (r) = (f(r) + f(-r)) /2$
is even and $f_o (r) = (f(r) - f(-r)) /2$ is odd. Evidently, we have
$\|f_e\|_{\kappa,p} \le \|f\|_{\kappa,p}$ and $\|f_o\|_{\kappa,p} \le 
\|f\|_{\kappa,p}$, we only need to prove the stated inequality for
$f$ being an odd function. 

Let $f$ be an odd function and define $g(r)  = f(r) /r$ for $r \ne 0$. 
Then $g$ is even.  The explicit formula of $\tau_r f$ in \eqref{eq:1.6} 
shows that 
$$
  \tau_r f(s) = (s-r) \tau_r g(s) = s \tau_r g(s) - r \tau_r g(s) 
$$
so that the Riesz operator can be written as a sum of two terms, 
$$
  R^\kappa f(s) = c_\kappa s \int_{\RR} \tau_r g(s) \frac{r} {|r|^2} dr -
    c_\kappa \int_{\RR} \tau_r g(s) dr : =       
   c_\kappa R_1^\kappa f(s) - c_\kappa R_2^\kappa f(s).
$$

For the first term we start with the following observation. Let 
$\kappa = (m-1)/2$, where $m \in \NN$. Define $F(x) = g(\|x\|)$ for 
$x \in \RR^m$. Let $\Omega$ be a harmonic polynomial of first degree 
on $\RR^m$ and let $K(x) = {\Omega(x)}/{\|x\|^{m+1}}$. We consider
the convolution $F* K$ in $L^1(\RR^m)$. Using the spherical-polar 
coordinates and the ordinary Funk-Hecke formula, we get
\begin{align*}
(F * K) (x) &= \int_0^\infty r^{m-1} \int_{S^{m-1}} 
g \left( \sqrt{\|x\|^2 + r^2 - 2 \|x\| \, r \langle x',y'\rangle}
 \right) \Omega(y') d\omega(y') \frac{dr}{r^m}\\ 
&  = c\, \Omega(x') \int_0^\infty \int_{-1}^1 g\left(
   \sqrt{\|x\|^2+ s^2 - 2\|x\| r t}\right)
   t (1-t^2)^{\frac{m-3}{2}} dt \frac{dr}{r} \\
&  = c\, \Omega(x') \int_{- \infty}^\infty \int_{-1}^1 g\left(
   \sqrt{\|x\|^2+ s^2 - 2\|x\| s t}\right)
   t (1-t^2)^{\frac{m-3}{2}} dt \sign(r) \frac{dr}{|r|}.
\end{align*}
Since $\sign (r) /|r| = r /r^2$ is odd in $r$, changing variables $t \to -t$ 
and $r \to -r$ shows that we can replace $t (1-t^2)^{m-3}/2$ by
$(1+t) (1-t^2)^{\kappa -1}$ in the last expression, so that the inner
integral becomes $\tau_r g (\|x\|)$, from which we get 
$$
 (F * K) (x) =
 c\, \Omega(x') \int_{-\infty}^\infty \tau_r f(\|x\|) \frac{r}{|r|^2} dr.  
$$
Hence, $\|x\| (F*K)(x)  =  c\, \Omega(x') R_1^\kappa f(\|x\|)$. Using the
spherical-polar coordinates and integrating over $\RR^m$ we get
\begin{align*}
 \int_{\RR^m}  |(F * K) (x)|^p \|x\|^p dx & =  
     c \int_0^\infty |R_1^\kappa f(s)|^p s^{m-1} ds 
      \int_{S^{d-1}} \Omega(x') d \omega(x') \\
   & =  c \int_0^\infty |R_1^\kappa f(s)|^p s^{2 \kappa} ds.
\end{align*} 
Since $R_1^\kappa f$ is an even function, this gives
\begin{align*}
  \int_{-\infty}^\infty |R_1^\kappa f(s)|^p |s|^{2 \kappa} ds  
 &  \le c  \int_{\RR^m}  |(F * K) (x)|^p \|x\|^p dx \\  
 &  \le c  \int_{\RR^m}  |F (x)|^p \|x\|^p dx 
\end{align*}
for $p > m /(m-1) = 2 \kappa +1 /(2\kappa)$ since $\|x\|^p \in A_p(\RR^m)$
for $p > m/(m-1)$. By definition, $F(x) = g(\|x\|)$ and $g(s) = f(s)/s$,  
so that 
$$
 \int_{\RR^m}  |F (x)|^p \|x\|^p dx = \int_{\RR^m} |f(\|x\|)|^p dx
   = \int_0^\infty |f(r)|^p r^{m-1} dr 
   = \frac{1}{2}\int_{- \infty}^\infty |f(r)|^p |r|^{2\kappa} dr.
$$
This takes care of the first term $R_1^\kappa f$ for 
$p > 2 \kappa +1/(2\kappa)$. 

For the second term we consider the operator $F \mapsto F * K_1$, where 
$F$ is as above and $K_1(x) = \|x\|^{1-m}$, which agrees with the notation 
in \eqref{eq:2.1}. Similar to $F*K$ we get  
\begin{align*}
  F*K_1(x) & = c \int_0^\infty \int_{-1}^1 g\left(
   \sqrt{\|x\|^2+ s^2 - 2\|x\| r t}\right)
    (1-t^2)^{\frac{m-3}{2}} dt dr \\
   & = c \int_{-\infty}^\infty \tau_r g(\|x\|) dr 
     = c R_2^\kappa(\|x\|).
\end{align*}   
Therefore, since $R_2^\kappa$ is even, 
\begin{align*}
 \int_{-\infty}^\infty |R_2^\kappa(s)|^p |s|^{2\kappa} ds & = 
  2 \int_0^\infty |R_2^\kappa(s)|^p s^{m-1} ds \\
 & = c \int_{\RR^m} |F*K_1 (x)|^p  dx   
  \le c \int_{\RR^m} |F (x)|^p \|x\|^p dx  
\end{align*}   
for $p > m / (m-1)$ using Proposition \ref{prop:2.2}. Again, the last 
integral is the same as $\|f\|_{\kappa,p}^p$, which takes care of the 
second term for $p > 2 \kappa+1/(2\kappa)$. Together, we have proved that 
$\|R^\kappa\|_{\kappa,p} \le   c \|f\|_{\kappa,p}$ or $p > 
2 \kappa+1/(2\kappa)$. For $1 < p \le (2\kappa+1)/(2\kappa)$ we use the
standard duality argument.
\end{proof}

\end{document}